\documentclass[submission]{dmtcs}

\usepackage[english]{babel}
\usepackage{amsmath}
\usepackage{amsfonts}

\newcommand\dedication{\section*{Dedication}}
\newtheorem{definition}{Definition}
\newtheorem{proposition}{Proposition}
\newtheorem{theoreme}{Theorem}
\newtheorem{lemme}{Lemma}
\newtheorem{corollaire}{Corollary}
\newtheorem{remark}{Remark}
\newtheorem{exemple}{Example}
\newtheorem{note}{Note}
\def\mref#1{(\ref{#1})}

\def\hs{\hbox to 3mm{}}
\def\hhs{\hbox to 5cm{}}
\def\ss{\smallskip}

\def\ncp#1#2{#1\langle #2 \rangle}
\def\ncs#1#2{#1\langle\langle #2 \rangle\rangle}

\def\Z{\mathbb{Z}}
\def\R{\mathbb{R}}
\def\C{\mathbb{C}}

\def\N{\mathbb{N}}
\def\K{\mathbb{K}}
\def\x{\mathtt{x}}
\def\RF{\mathsf{RF}}
\def\End{\mathsf{End}}
\def\Id{\mathsf{Id}}

\def\End{\mathsf{End}}
\def\binomial#1#2{\left(\,\begin{matrix}#1 \\ #2\end{matrix}\,\right)}
\def\ra{\rightarrow}

\def\al{\alpha}

\def\la{\lambda}
\def\Om{\Omega}

\author{G\'erard H. E. Duchamp \addressmark{1}\thanks{gerard.duchamp@lipn.univ-paris13.fr}, Laurent Poinsot \addressmark{1}\thanks{laurent.poinsot@lipn.univ-paris13.fr}, Allan I. Solomon \addressmark{4,2}, Karol A. Penson \addressmark{2},
Pawel Blasiak, Andrzej Horzela \addressmark{3}.}
\title{Ladder Operators and Endomorphisms in Combinatorial Physics}
\address{\addressmark{1}LIPN, UMR CNRS 7030.
Institut Galil\'ee - Universit\'e Paris-Nord,\\
99, avenue Jean-Baptiste Cl\'ement
93430 Villetaneuse, France.\\
\addressmark{2} Universit\'{e} Pierre et Marie Curie,
Laboratoire  de  Physique   Th\'{e}orique  de la Mati\`ere Condens\'ee, CNRS UMR 7600\\
Tour 16, $2^{i\grave{e}me}$ \'{e}tage, 4, place Jussieu, F 75252 Paris Cedex 05,
France.\\
\addressmark{3} H. Niewodnicza{\'n}ski Institute of Nuclear Physics,
Polish Academy of Sciences\\
Division of Theoretical Physics\\
ul. Radzikowskiego 152, PL 31-342 Krak{\'o}w, Poland.\\
\addressmark{4}
The Open University\\
 Physics and Astronomy Department\\
Milton Keynes MK7 6AA, United Kingdom
}
\revision{26}
\received{\today}
\revised{\today}
\begin{document}
\maketitle
\date{}
\begin{abstract}
Starting with the Heisenberg-Weyl algebra, fundamental to quantum physics, we first show how the ordering of the non-commuting operators intrinsic to that algebra gives rise to generalizations of the classical Stirling Numbers of Combinatorics. These may be expressed in terms of infinite, but {\em row-finite}, matrices, which may also be considered as endomorphisms of $\C[[\mathtt{x}]]$.  This  leads us to consider endomorphisms in more general spaces, and these in turn may be expressed in terms of generalizations of the ladder-operators familiar in physics.
\end{abstract}
\tableofcontents

\section{Introduction}

Combinatorics has a long history which can be traced back to the times when Greek, Chinese and Persian mathematicians (to name but a few) began with this particular and fruitful blend of configuration and counting.

More recently, due to the great masters of the past (Euler, Bernouilli etc ...), this ``art of counting'' avoided the fate of becoming a ``collections of recipes'' and, under the impetus of the modern fields of Algorithms and Computer Sciences, acquired its {\em Letters Patent}  and so pervaded many domains of Classical Sciences, such as Mathematics and Physics.

In return, the sciences which interact with Combinatorics can transmit to the latter some of their art.\\
This is the case of the emerging field of  ``Combinatorial Physics'' which has the potential of revitalizing  mathematical   features that have been  familiar to physicists for over a century, such as  tensor calculus, structure constants, operator calculus, infinite matrices, and so on.

%The goal of this paper is to concentrate on two aspects of this interaction: namely, on the one hand, the use of %one-parameter groups  to solve counting problems and, on the other hand,  series and infinite matrices to represent %them.

In this paper we describe one aspect of this interaction; namely, how well-known concepts in quantum physics such as creation and annihilation operators, and ladder operators, translate to combinatorial "counting" ideas as exemplified by Stirling numbers, which may find their expression in terms of infinite matrices.  Such an infinite matrix is more  generally to be thought of as a (linear) transformation from a  linear space to itself, that is, a linear   endomorphism. This is also a rigorous context in which to describe the traditional {\em ladder operators} of physics.

The paper is organized as follows: We start by introducing the well-known Heisenberg-Weyl associative algebra generated by the creation and annihilation operators of second-quantized  physics; this is   a graded algebra.  Consideration of exponentials of elements of this algebra   leads one to a generalization of the classical combinatorial Stirling numbers, as well as one-parameter groups - crucial in quantum physics. Arrays of such numbers lead us to the algebra of row-finite infinite matrices.  We then consider linear endomorphisms as a natural sequel to these matrices, and their representations.  We relate these to a generalization of the idea of ladder operators, and conclude by giving some results concerning the relation between endomorphisms and ladder operators.

\acknowledgements The authors wish to acknowledge support from Agence Nationale de la Recherche (Paris, France) under Program No. ANR-08-BLAN-0243-2 and from PAN/CNRS Project PICS No.4339(2008-2010). Two of us (P.B. and A.H.) wish to acknowledge support from Polish Ministry of Science and Higher Education under Grants Nos. N202 061434 and N202 107 32/2832. Two others (G.H.E.D. and L.P.) would like to acknowledge support from ``Projet interne au LIPN 2009'' ``Polyz\^eta functions''.

\dedication

{\em Philippe Flajolet has not only blazed a successful trail in his own areas of Combinatorics but, through  his
support given generously to others,  has  stimulated new branches to develop.  On the occasion of
Philippe's  60th birthday we are happy to acknowledge with gratitude his encouragement to the
development of our own field of Combinatorial Physics, and look forward to many future years of
professional engagement.}

\section{The Heisenberg-Weyl algebra}
\subsection{Formal definition}

In Quantum Physics \cite{BJ,BHJ,WH} and more recently in Combinatorics \cite{Fomin,Stanley2} and Combinatorial Physics \cite{GOF13}, one often encounters pairs of operators $(A,B)$ such that
\begin{equation}\label{hw-relation}
	AB-BA=I
\end{equation}
 where $I$ stands for the identity in some associative algebra. The appearance of this relation in 1925 at once forced Born, Heisenberg and Jordan to the consideration of infinite matrices. Indeed, it can be shown that relation \mref{hw-relation} cannot be represented by (finite) matrices with elements in a field of characteristic zero (simply take the trace of each side). The first choice of faithful representation for \mref{hw-relation} is with (densely defined) unbounded operators in a Hilbert space (traditional Fock space) or with continuous operators in a Fr\'echet space \cite{Bou-TVS,GOF4,Tr}.\\
 One can  formally define the Heisenberg-Weyl algebra by

\begin{equation}
	HW_\C=\ncp{\C}{b,b^+}/\mathcal{J}_{HW}
\end{equation}
where $\ncp{\C}{b,b^+}=\C\hspace{1mm} [\{b,b^+\}^*]$ is the algebra of the free monoid $\{b,b^+\}^*$ \cite{BP,BR,Lothaire} i. e. the algebra of non-commutative polynomials; and $\mathcal{J}_{HW}$ is the two-sided ideal generated by $(bb^+-b^+b-1)$. Note that this definition, together with the arrow
\begin{equation}\label{canarrow}
	s: \ncp{\C}{b,b^+}\ra HW_\C
\end{equation}
clears up all the ambiguities concerning normal forms (Normal ordering \cite{SBDHP} and the so-called ``double dot'' operation) which are traditional in Quantum Physics. From now on, we set $a=s(b)$ and $a^+=s(b^+)$.\\
In general, by the \emph{normal ordering} \cite{GOF12} of a general expression $F(a^{\dag},a)$ we mean
$F^{(n)}(a^{\dag},a)$ which is obtained by moving all the annihilation operators $a$ to the right \emph{using} the
commutation relation of Eq.(\ref{hw-relation}). This procedure yields an operator whose action is equivalent to the original one, \emph{i.e.} $F^{(n)}(a^{\dag},a)=F(a^\dag,a)$ as operators, although the form (which lives in $\ncs{\C}{b,b^+}$ or $\ncp{\C}{b,b^+}$ ) of the expressions in terms of $a$ and $a^\dag$ may be completely different.
From \mref{canarrow}, it is an easy exercise to prove that $\Big((a^+)^ia^j\Big)_{i,j\in\N}$ is a basis of $HW_\C$ (basis of the normal order).\\
On the other hand the \emph{double dot} operation $:\!F(a^{\dag},a)\!:$
consists of
applying the same ordering procedure but \emph{without} taking into account the
commutation relation of Eq.\mref{hw-relation}, {\it i.e.} moving all annihilation
operators $a$
to the right as if they commuted with the creation operators $a^\dag$. The structure constants are given in \cite{GOF13} and can be obtained from the following formula \footnote{This formula can also easily be derived from the ``rook'' equivalent of Wick's theorem \cite{Varvak}. Note that the summation index $k$ ranges in the interval $[0..min(j_1,i_2)]$.}
\begin{equation}
	(a^+)^{i_1}a^{j_1}(a^+)^{i_2}a^{j_2}=\sum_{k\geq 0}k!\binomial{j_1}{k}\binomial{j_2}{k} (a^+)^{i_1+i_2-k}a^{j_1+j_2-k}\ .
\end{equation}

\subsection{Grading of the Heisenberg-Weyl algebra}

Setting, for $e\in \Z$
\begin{equation}
	HW_\C^{(e)}=span_\C((a^+)^ia^j)_{i-j=e}
\end{equation}
one has
\begin{equation}
	HW_\C=\bigoplus_{e\in \Z} HW_\C^{(e)} \textrm{ and } HW_\C^{(e_1)}HW_\C^{(e_2)}\subset HW_\C^{(e_1+e_2)}
\end{equation}
for all $e_1,e_2\in \Z$. This natural grading makes $HW_\C$ a $\Z$-graded algebra. One often uses the following (faithful) representation $\rho_{BF}$ by operators on $\C[[\mathtt{x}]]$.

\begin{eqnarray}\label{Bargmann-Fock}
\left\{\begin{array}{l}
\rho_{BF}(a)=\frac{d}{d\mathtt{x}}\\\\
\rho_{BF}(a^+)=(S\mapsto \mathtt{x}S)\ .
\end{array}\right.
\end{eqnarray}
This representation, known as the Bargmann-Fock representation is graded for the preceding grading as, when restricted to $\C[\mathtt{x}]$, $\rho_{BF}(a)$ is of degree $-1$ and $\rho_{BF}(a^+)$ of degree $1$.

\ss
In general, and more concretely, we may associate many important operators of quantum physics with elements of $HW_\C$. In particular,  an element $\Omega\in HW_\C$ being given, one would like to consider the evolution group \cite{D1}
$$
\Big(e^{\la\Om}\Big)_{\la \in \R}\ .
$$
For example, such one-parameter groups are important in quantum dynamics, where the parameter $\la$ is the time $t$; or in quantum statistical mechanics, where $\la$ is the negative inverse temperature.

Some questions which arise are \\ \\
Q1) Is this group well defined ? through which representation ? what is the domain ?\\
Q2) Which combinatorial methods may be extracted from knowledge of this group ?

\bigskip
Our first task  is to get the normal order of the powers $\Om^n$.

\section{Combinatorics of infinite matrices}
\subsection{Homogeneous operators and generalized Stirling numbers}

Before defining representations (or realizations) of the one-parameter group $\Big(e^{\la\Om}\Big)_{\la \in \R}$, one can consider the problem of normal ordering the powers of $\Om$
\begin{equation}
	\mathcal{N}(\Om^n)=\sum_{i,j\in \N} \al(n,i,j) (a^+)^ia^j\ .
\end{equation}
In general this is  a three-parameter problem but, taking advantage of the preceding gradation, one can start with a homogeneous operator of degree (or excess) $e$
\begin{equation}\label{eq9}
	\Om=\sum_{i-j=e} \al(i,j) (a^+)^ia^j
\end{equation}
and remark that

\begin{eqnarray}\label{Bargmann-Fock2}
\mathcal{N}(\Om^n)=\left\{\begin{array}{l}
(a^+)^{ne}\sum_{k=0}^\infty S_{\Om}(n,k) (a^+)^ka^k\  \textrm{; if $e\geq 0$}
\\\\
\Big(\sum_{k=0}^\infty S_{\Om}(n,k) (a^+)^ka^k\Big)a^{n|e|}\  \textrm{; if $e<0$}
\end{array}\right.
\end{eqnarray}

which was used  as the definition of  ``Generalized Stirling Numbers'' as introduced in \cite{BPS1,BPS2} for strings and generalized to  homogeneous operators in \cite{GOF4} (see also \cite{MBP}). These numbers recently attracted the attention of Combinatorialists \cite{Flajolet} who found it a nontrivial generalization of numbers known for some 200 years \cite{SloaneStirling}.\\
The reason for the  name {\em Stirling Numbers} lies in the first example below, following which  we give two more examples.

\smallskip
For $\Omega=a^+a$, one gets the usual matrix of Stirling numbers of the second kind.
\begin{equation}
\left\lceil
{\begin{array}{rrrrrrrr}
1 & 0 & 0 & 0 & 0 & 0 & 0 &\cdots\\
0 & 1 & 0 & 0 & 0 & 0 & 0 &\cdots\\
0 & 1 & 1 & 0 & 0 & 0 & 0 &\cdots\\
0 & 1 & 3 & 1 & 0 & 0 & 0 &\cdots\\
0 & 1 & 7 & 6 & 1 & 0 & 0 &\cdots\\
0 & 1 & 15 & 25 & 10 & 1 & 0 &\cdots\\
0 & 1 & 31 & 90 & 65 & 15 & 1&\cdots\\
\vdots & \vdots & \vdots  & \vdots  & \vdots  & \vdots  & \vdots &\ddots\\
\end{array}}
 \right.
\end{equation}

\smallskip
For $\Omega=a^+aa^+$, we have
\begin{equation}
\left\lceil
{\begin{array}{rrrrrrrr}
1 & 0 & 0 & 0 & 0 & 0 & 0 &\cdots\\
1 & 1 & 0 & 0 & 0 & 0 & 0 &\cdots\\
2 & 4 & 1 & 0 & 0 & 0 & 0 &\cdots\\
6 & 18 & 9 & 1 & 0 & 0 & 0 &\cdots\\
24 & 96 & 72 & 16 & 1 & 0 & 0 &\cdots\\
120 & 600 & 600 & 200 & 25 & 1 & 0 &\cdots\\
720 & 4320 & 5400 & 2400 & 450 & 36 & 1&\cdots\\
\vdots & \vdots & \vdots  & \vdots  & \vdots  & \vdots  & \vdots &\ddots\\
\end{array}}
\right.
\end{equation}

\smallskip
For $\Omega=a^+aaa^+a^+$, one gets
\begin{equation}
\left\lceil
{\begin{array}{rrrrrrrrrr}
1 & 0 & 0 & 0 & 0 & 0 & 0 & 0 & 0 & \cdots\\
2 & 4 & 1 & 0 & 0 & 0 & 0 & 0 & 0 &\cdots\\
12 & 60 & 54 & 14 & 1 & 0 & 0 & 0 & 0 &\cdots\\
144 & 1296 & 2232 & 1296 & 306 & 30 & 1 & 0 & 0 &\cdots\\
2880 & 40320 & 109440 & 105120 & 45000 & 9504 & 1016 & 52 & 1 &\cdots\\
\vdots & \vdots & \vdots  & \vdots  & \vdots  & \vdots  & \vdots & \vdots & \vdots &\ddots\\
\end{array}}
\right.
\end{equation}

In any case, the matrix $(S(n,k))_{n,k\in\N}$ has all its rows $(S(n,k))_{k\in\N}$ finitely supported. We call these matrices ``row-finite'' \cite{GOF4,Pi}.

\ss
We will see in the next paragraph that the ``row-finite'' matrices form a very important algebra which we denote by $\RF(\N;\C)$ in the sequel.

\subsection{An excursion to topology: transformation of sequences}

Let $\mathbf{d}=(d_n)_{n\in\N}$ be a set of non-zero complex denominators. To each row-finite matrix $(M[n,k])_{n,k\in\N}$, one can associate an operator $\Phi_M\in\End(\C[[\mathtt{x}]])$ such that the image of $f=\sum_{k\in\N}a_k\frac{\mathtt{x}^k}{d_k}\in\C[[\mathtt{x}]]$ is defined by

\begin{equation}
	\Phi_M(f)=\sum_{n\in\N}b_n\frac{\mathtt{x}^n}{d_n} \textrm{; with } b_n=\sum_{k\in\N}M[n,k]a_k\ .
\end{equation}
Note that if we endow $\C[[\mathtt{x}]]$ with the Fr\'echet topology of simple convergence of the coefficients (this structure is sometimes called the ``Treves topology'', see \cite{Tr})     \emph{i.e.}, defined by the seminorms
\begin{equation}
	p_n(f):=|a_n| \textrm{; with } f=\sum_{k\in\N}a_k \mathtt{x}^k
\end{equation}
with each  $\Phi_M$  continuous; then the following proposition states that there is no other case:

\begin{proposition} The correspondence $M\ra \Phi_M$ from $\RF(\N;\C)$ to $\mathcal{L}(\C[[\mathtt{x}]])$ (continuous endomorphisms) is one-to-one and linear. Moreover $\Phi_{MN}=\Phi_M\circ \Phi_N$.
\end{proposition}
\begin{proof}
The proof of this proposition is not difficult and left to the reader.
\end{proof}

\ss
As an application of the preceding, one can remark that, through the Bargmann-Fock representation $\rho_{BF}$, the one parameter group $e^{\lambda\Omega}$ always makes sense for homogeneous operators (as defined in Eq. \mref{eq9}) since the matrix $\phi^{-1}(\rho_{BF}(\Omega))$\footnote{These matrices are different from the ``Generalized Stirling matrices'' defined by Eq. \mref{Bargmann-Fock2}. Their non-zéro elements are supported by a line parallel to the diagonal.} is

\begin{itemize}
	\item strictly upper-triangular when $e<0$
	\item diagonal when $e=0$
	\item strictly lower-triangular when $e>0$ .
\end{itemize}

Then $e^{\lambda\Omega}$ is meaningful as (a group of) operators on appropriate spaces.

\subsection{One-parameter groups and Stirling matrices}

In this paragraph we focus on the combinatorics of operators containing at most one annihilator (in this context $d/dx$) so that $\rho_{BF}(\Omega)$ is of the form
\begin{equation}\label{affine field}
	q(x)\frac{d}{dx} + v(x)
\end{equation}
(sum of a scalar field and a true vector field). One-parameter groups generated by these operators can, of course, be integrated using PDE \cite{D1} but, here we give a ``congugacy trick'' which aims at proving that an operator of the type \mref{affine field} is conjugated to the vector field $q(x)\frac{d}{dx}$.\\
So, to compute $e^{\lambda (q(x)\frac{d}{dx} + v(x))}[f]$, one can use the following procedure ($q$ and $v$ are supposed to be at least continuous). We first take $v\equiv 0$ (vector field case)

\begin{itemize}
	\item if $q\equiv 0$ (and $v\equiv 0$) then $e^{\lambda \rho_{BF}(\Omega)}[f]=f$ (trivial action) ;
	\item if $q\not\equiv 0$ then choose an open interval $I\neq\emptyset$ in which $q$ never vanishes and $x_0\in I$ ;
	\item for $x\in I$ set
\begin{equation}
	F(x)=\int_{x_0}^x \frac{dt}{q(t)}
\end{equation}
and set $J=F(I)$ (open interval). Then $F:I\ra J$ is one-to-one (as $F$ is strictly monotonic) ;
	\item for suitable $(x,\lambda)$, set
\begin{equation}
	s_\lambda(x)=F^{-1}(F(x)+\lambda).
\end{equation}
$s_\lambda$ is a deformation of the identity since $(x,\lambda)\mapsto s_\lambda(x)$ is continuous (and even of class $C^1$) on its domain and $s_0(x)=x$ ;
	\item for small values of $\lambda$ , $e^{\lambda (q(x)\frac{d}{dx})}$ coincides with the substitution $f\mapsto f\circ s_\lambda$. To see this, it is sufficient to remark that the exponential of a derivation (such as $\lambda (q(x)\frac{d}{dx})$) is an automorphism, which means a substitution in the (test) function spaces under consideration.
\end{itemize}

Now, one can indicate how to integrate the one-parameter group $e^{\lambda (q(x)\frac{d}{dx} + v(x))}$ for general $v$.
($I,\ F,\ s_\lambda$ are as above).

\begin{itemize}
	\item On $I$, set
\begin{equation}
	u(x)=e^{\int_{x_0}^x \frac{v(t)}{q(t)}dt}\ ;
\end{equation}
	\item one checks easily that
\begin{equation}
	\rho_{BF}(\Omega)=(q(x)\frac{d}{dx} + v(x))=\frac{1}{u}(q(x)\frac{d}{dx})u
\end{equation}
in the sense that, on each function in its domain $\rho_{BF}(\Omega)$ operates as the composition of
\begin{itemize}
	\item multiplication of $f$ by $u$
	\item action of the vector field $(q(x)\frac{d}{dx})$ (now on $uf$)
	\item division by $u$ ;
\end{itemize}
	\item then, using the fact that exponentiation commutes with conjugacy, the exponential reads
	\begin{equation}
	e^{\lambda(q(x)\frac{d}{dx} + v(x))}={u^{-1}}e^{\lambda(q(x)\frac{d}{dx})}u\ .
	\end{equation}
\end{itemize}

Using the preceding definitions, the action now takes the form
\begin{equation}
U_\lambda[f](x)=e^{\lambda(q(x)\frac{d}{dx} + v(x))}[f](x)=\frac{u(s_\lambda(x))}{u(x)} f(s_\lambda(x))\ .
\end{equation}

One can check \textit{a posteriori} the validity of this procedure, using a tangent vector technique as follows
\begin{itemize}
	\item check that, for small values of $\lambda,\theta$, one has
\begin{equation}
U_\lambda\circ U_\theta=U_{\lambda+\theta}\ ;	
\end{equation}
	\item check that
\begin{equation}
	\frac{d}{d\lambda}\Big|_{\lambda=0}U_\lambda[f](x)=(q(x)\frac{d}{dx}+v(x))f(x)\ .
\end{equation}
\end{itemize}

\begin{remark} Transformations of type
\begin{equation}
	f\ra g.(f\circ s)
\end{equation}
are called  {\em substitutions with prefunctions } in combinatorial physics\cite{GOF4}. It can be shown that under nice conditions ($g, s$ analytic in a neighbourhood of the origin, $g(0)=1, s=x+\cdots$) these transformations form a (compositional) Lie group (infinite dimensional of Fr\'echet type, see \cite{GOF4}). The infinitesimal generators of these transformations are precisely of the form $q(x)\frac{d}{dx} + v(x)$.
\end{remark}

We now give an example of integration of the one-parameter group $e^{\lambda\rho_{BF}(\Omega)}$ for
\begin{equation}
\Omega=(a^+)^2aa^++a^+a(a^+)^2.	
\end{equation}

\begin{exemple}\label{ex1}
One has the conjugated form
\begin{equation}
\rho_{BF}(\Omega)=x^2\frac{d}{dx}x+	x\frac{d}{dx}x^2=x^{-\frac{3}{2}}(2x^3\frac{d}{dx})x^{\frac{3}{2}}\ .
\end{equation}
Using the procedure described above, one obtains the one-parameter group of transformations $U_\lambda$
\begin{equation}
	U_\lambda [f](x)=\sqrt[4]{\frac{1}{(1-4\lambda x^2)^3}}\times f(\sqrt{\frac{x^2}{1-4\lambda x^2}})\ .
\end{equation}
The reader is invited to check that, for suitably small values of the parameters\\ (i.e. $|\lambda|+|\theta|<\frac{1}{4x^2}\leq +\infty$), $U_\lambda\circ U_\theta=U_{\lambda+\theta}$ by direct computation.
\end{exemple}

Once integrated, the one-parameter group $U_\lambda$ reveals the Generalized Stirling matrix as expressed by the following result.
\begin{proposition} With the definitions introduced and $e\geq 0$, the two following conditions are equivalent\\
(where $f \ra U_\lambda[f]$ is the one-parameter group $exp(\lambda \rho_{BF}(\Omega))$.\\
i)
\begin{equation}
	\sum_{n,k\geq 0} S_\Omega(n,k)\frac{x^n}{n!}y^k=g(x)e^{y\phi(x)}
\end{equation}
ii)
\begin{equation}\label{second_cond}
	U_\lambda[f](x)=g(\lambda x^e)f(x(1+\phi(\lambda x^e)))
\end{equation}
\end{proposition}

\begin{proof} One first has the following equality between continuous operators
\begin{eqnarray}\label{prepa_calc}
U_\lambda=\sum_{n,k\geq 0} S_\Omega(n,k)\frac{\lambda^n}{n!} x^{ne}  x^k(\frac{d}{dx})^k .
\end{eqnarray}
%Now, as the two members of \mref{second_cond} are continuous and linear in $f$, it is enough to test it on a total subset \cite{Bou-TVS} of functions for the Treves topology\footnote{The usual - ultrametric - topology would not be enough for $e=0$.} for example the monomials (here ) . 
Assuming (i), let us check (ii) for $f$ a monomial (i. e. choose the test functions $f=x^j$, for $j=0,1,\cdots$) 
\begin{eqnarray}\label{point_d'appui}
U_\lambda (x^j)&=&\sum_{n\geq 0}\sum_{k=0}^j S_\Omega(n,k)\frac{(\lambda x^e)^n}{n!}\frac{j!}{(j-k)!}x^j=\cr
               &=& x^j \sum_{k=0}^j \Big( [y^k] g(\lambda x^e)e^{y\phi(\lambda x^e)}\Big)\frac{j!}{(j-k)!}=\cr
               &=& g(\lambda x^e) x^j \sum_{k=0}^j \binomial{k}{j} \phi(\lambda x^e)^k=
               g(\lambda x^e)\Big(x\big(1+\phi(\lambda x^e)\big)\Big)^j\ .
\end{eqnarray}
Now as the two members of \mref{second_cond} are continuous and linear in $f$ and the set of monomials is total \cite{Bou-TVS} in the space of formal power series endowed with the Treves topology\footnote{The usual - ultrametric - topology would not be enough for $e=0$.}, we have (ii).\\
Conversely, if one assumes (ii), one has 
\begin{equation}
	U_\lambda (e^{yx})=g(\lambda x^e)e^{yx(1+\phi(\lambda x^e))}
\end{equation}
 and, from \mref{point_d'appui}, one gets
\begin{equation}
\sum_{n,k\geq 0} S_\Omega(n,k)\frac{(\lambda x^e)^n}{n!} (xy)^k=g(\lambda x^e)e^{yx\phi(\lambda x^e)}.	
\end{equation}
A legitimate change of variables ($\lambda x^e\ra x\ ;\ xy\ra y$) gives (i).
\end{proof}

{\bf Example \ref{ex1} continued}\\
{\it With $\Omega=(a^+)^2aa^++a^+a(a^+)^2$, one has the one-parameter group
\begin{equation}
	U_\lambda[f](x)=\sqrt[4]{\frac{1}{(1-4\lambda x^2)^3}}\times f(\sqrt{\frac{x^2}{1-4\lambda x^2}})\ .
\end{equation}
Then, applying the preceding correspondence, one gets
\begin{equation}
	\sum_{n,k\geq 0} S_\Omega(n,k)\frac{x^n}{n!}y^k=\sqrt[4]{\frac{1}{(1-4x)^3}}\ \ e^{y(\sqrt{\frac{1}{(1-4x)}}-1)}=
	\sqrt[4]{\frac{1}{(1-4x)^3}}\ \ e^{y(\sum_{n\geq 1} c_nx^n)}
\end{equation}
where $c_n=\binomial{2n}{n}$ are the central binomial coefficients.
}

\section{Representation of endomorphisms in more general spaces}

\subsection{Notation}

Consider $\K$ a (commutative) field and $\K[\x]$ the $\K$-vector space of polynomials in the indeterminate
$\x$. Denote by $\End(V)$ the algebra of linear endomorphisms of any $\K$-vector space $V$. If $\phi$ and $\psi$ are both elements of $\End(V)$, then with $\phi\psi$ denoting the usual composition ``$\phi\circ\psi$" of linear mappings, we have for any integer $n$
\begin{equation}
\phi^n := \left \{
\begin{array}{lll}
\Id_V & \mbox{if} & n=0\ ,\\
\underbrace{\phi\circ \cdots \circ \phi}_{n\ \mathit{times}} & \mbox{if} & n>0
\end{array}
\right .
\end{equation}
where $\Id_V$ is the identity mapping of $V$.
Let $\mathbf{e}:=(e_i)_{i\in I}$ be a basis of $V$ ($V$ which we assume does  not reduce to $(0)$). We denote the decomposition of any vector $v \in V$ with respect to $\mathbf{e}$ by
\begin{equation}
\displaystyle\sum_{i\in I}\langle v,e_i\rangle e_i
\end{equation}
where\footnote{The notation ``$\langle v,w \rangle$" is commonly referred to as a ``Dirac bracket". It was successfully used  (for the same reason of duality) by Sch\"utzenberger to develop his theory of automata \cite{BP,BR,Eil74}.} $\langle v,e_i\rangle$ is the coefficient of the projection of $v$ onto the subspace $\K e_i$ generated by $e_i$ in $V$. Obviously, all but a finite number of the coefficients $\langle v,e_i\rangle$ are equal to zero. If $(I,\leq)$ is a linearly ordered (nonempty) set bounded from below (with $\widehat{0}$ as its minimum\footnote{We follow the notation of \cite{Stanley} for the lowest element.}), and, if $v\not=0$, then the {\emph{degree}} of $v$ (with respect to $\mathbf{e}$) is defined by
\begin{equation}
\deg_{\bf{e}}(v) := \max\{i\in I : \langle v,e_i\rangle\not=0\}
\end{equation}
and
\begin{equation}
\deg_{\bf{e}}(0):=-\infty\
\end{equation}
where $-\infty \not\in I$, and the relation $-\infty < i$ for each $i\in I$ extends the order of $I$ to $\overline{I}:=I\cup\{-\infty\}$. If $v\not=0$, then the nonempty finite set $\{i\in I : \langle v,e_i\rangle\not=0\}$ admits a greatest element, since $I$ is totally ordered, so that $\deg_{\bf{e}}(v)$ is well-defined.
Thus, the following equality holds (for any $v\not=0$)
\begin{equation}
v=\displaystyle\sum_{\widehat{0}\leq i\leq \deg_{\mathbf{e}}(v)}\langle v,e_i\rangle e_i
\end{equation}
with $\langle v,e_{\deg_{\mathbf{e}}(v)}\rangle \not=0$. In particular, taking $\mathbf{x}:=(\x^n)_{n\geq 0}$ as a basis of $\K[\x]$, any nonzero polynomial $P$ may be  written as the sum
\begin{equation}
P=\displaystyle\sum_{n=0}^{\deg(P)}\langle P,\x^n\rangle \x^n
\end{equation}
where  $\deg(P)$ is the usual degree  of $P$.

\subsection{Review of the classical result}\label{sec3}

It has been known since the paper of Pincherle and Amaldi \cite{PA01} that, for  a field $\K$ of characteristic zero, any linear endomorphism $\phi \in \End(\K[\x])$ may be expressed as the sum of a converging series in the operator $X$ of multiplication by the variable $\x$ and in the (formal) derivative (of polynomials) $D$. In \cite{KM86} (see also \cite{DiBL96} for some generalizations)  Kurbanov and Maksimov give  an explicit formula - recalled below -  for this sum.
\begin{theoreme}[\cite{KM86}]
Suppose that $\K$ is a field of characteristic zero. Let $\phi\in \End(\K[\x])$. Then $\phi$ is the sum of the summable series (in the topology of simple convergence on $\End(\K[\x])$ with $\K[\x]$ discrete) $\displaystyle\sum_{k=0}^{+\infty} P_k(X)D^k$ where $(P_k(\x))_{k\in\N}$ is a sequence of polynomials which satisfies the following recursion equation:
\begin{equation}
\begin{array}{lll}
P_0(\x)&=&\phi (1)\ ,\\
P_{n+1}(\x)&=&\displaystyle\phi(\frac{\x^{n+1}}{(n+1)!})-\sum_{k=0}^n P_k(\x)\frac{\x^{n+1-k}}{(n+1-k)!}\ .
\end{array}
\end{equation}
\end{theoreme}

In what follows, we generalize this result to any $\K$-vector space with a countable basis using a pair of rather general ladder operators instead of the usual ones, namely $X$ and $D$. The basic  idea is   to use only those operator properties which make possible an expansion similar to the classical case.

\subsection{Endomorphism expansion in terms of ladder operators}

From now on,  except for Example~\ref{evaluation_exmpl}, the field $\K$ is not assumed to be of characteristic zero. Let us consider a $\K$-vector space $V$ of  countable dimension. Let $\mathbf{e}:=(e_n)_{n\in\N}$ be an algebraic basis for this space. We can define two kinds of {\emph{ladder operators}} with respect to $\mathbf{e}$, namely, a {\emph{lowering operator}} $L_{\bf{e}}\in\End(V)$, by
\begin{equation}
\left \{
\begin{array}{lll}
L_{\bf{e}}e_0 &=& 0\ ,\\
L_{\bf{e}}e_{n+1}&=&e_{n}
\end{array}\right .
\end{equation}
and, a {\emph{raising operator}} $R_{\bf{e}}\in\End(V)$, by
\begin{equation}
R_{\bf{e}} e_n = e_{n+1}\ .
\end{equation}
Such operators were discussed  by Katriel and Duchamp \cite{KD95} as well as
Dubin, Hennings and Solomon \cite{DubHS1,DubHS2} in a more general context, and are  similar to the  creation and annihilation operators acting on an interacting Fock space of Accardi and Bo\.{z}ejko \cite{AB98}.
The operators $L_{\bf{e}}$ and $R_{\bf{e}}$ may also be regarded  as  the operators $D$ and $U$ described by Fomin in \cite{Fomin}, associated with the oriented graded graph $e_0 \leftarrow e_1 \leftarrow e_2 \leftarrow \cdots$ and $e_0 \rightarrow e_1 \rightarrow e_2 \rightarrow \cdots$.
\begin{definition}
Let $P\in\K[\x]$ and $\mathbf{u}:=(u_n)_{n\in\N}$ be a sequence of elements of $V$. We define $P(\mathbf{u})\in V$ by
\begin{equation}
P(\mathbf{u}):=\displaystyle\sum_{n\geq 0}\langle P,\x^n\rangle u_n = \sum_{n=0}^{\deg(P)}\langle P,\x^n\rangle u_n\ .
\end{equation}
\end{definition}
\begin{lemme}\label{isomorphisme}
Let $\mathbf{e}=(e_n)_{n\in \N}$ be a basis of $V$. The mapping
\begin{equation}
\begin{array}{llll}
\Phi_{\mathbf{e}} : & \K[\x] & \rightarrow & V\\
& P & \mapsto & P(\mathbf{e})
\end{array}
\end{equation}
is a linear isomorphism.
\end{lemme}
\begin{proof}
Straightforward.
\end{proof}
\begin{lemme}\label{petitlem}
Let $\mathbf{e}=(e_n)_{n\in\N}$ be a basis of $V$ and $R_{\mathbf{e}}$ be the raising operator associated with $\mathbf{e}$. For any polynomial $P\in\K[\x]$ we can define the operator $P(R_{\mathbf{e}}):=\displaystyle\sum_{n\geq 0}\langle P,\x^n\rangle R_{\mathbf{e}}^n$. Then we have
\begin{equation}
P(R_{\mathbf{e}})e_{0}=P(\mathbf{e})\ ,
\end{equation}
thus 
\begin{equation}
R_{\mathbf{e}}^n e_0=e_n\ .
\end{equation}
\end{lemme}
\begin{proof}
Omitted.
\end{proof}

Now suppose that $V$ is discrete (as is $\K$) and $\End(V)$, as a subspace of $V^V$,  is endowed with  the topology of compact convergence; that is, in this case, the topology of simple convergence (since the compact subsets of discrete $V$  are its finite subsets). As a result, $\End(V)$ becomes a complete topological $\K$-vector space (and even  a complete topological  $\K$-algebra). Using this topology we may consider summable families of operators on $V$. \\

We recall here some basics about summability in a general setting. Let $G$ be a Hausdorff commutative group, $(g_i)_{i\in I}$ a family of elements of $G$. An element $g\in G$ is the {\emph{sum}} of
$(g_i)_{i\in I}$ if, and only if, for each neighbourhood $W$ of $g$ there exists a finite subset $J_W$ of $I$ such that
\begin{equation}
\displaystyle\sum_{j\in J}g_j \in W
\end{equation}
for every finite subset $J\subset I$ containing $J_W$. The sum $g$ of a summable family $(g_i)_{i\in I}$ of elements of $G$ is usually denoted by
\begin{equation}
\displaystyle\sum_{i\in I}g_i\ .
\end{equation}
It is well-known that if $(g_i)_{i\in I}$ is a summable family with sum $g$, then for any permutation $\sigma$ of $I$, $g$ is also the sum of
$(g_{\sigma(i)})_{i\in I}$. When $G$ is complete, the following  condition (Cauchy) is equivalent to summability. A family $(g_i)_{i\in I}$ of $G$ satisfies {\emph{Cauchy's condition}} if, and only if, for every neighbourhood $W$ of zero there is a finite subset $J_W$ of $I$ such that
\begin{equation}
\displaystyle\sum_{k\in K}g_k\in W
\end{equation}
for every finite subset $K$ of $A$ disjoint from $J_W$. Many other properties and results about summable families may be found in \cite{Bou-GT}.\\

For instance, let $\mathbf{e}=(e_n)_{n\in \N}$ be a basis of $V$. Then for any sequence
$(\phi_n)_{n=0}^{\infty}\in \End(\K[\x])^{\N}$ of elements of $\End(V)$, the family $(\phi_n L^n_{\mathbf{e}})_{n\in\N}$ is easily shown to be summable. Due to the choice of  topology, the fact that $\mathbf{e}$ is a basis of $V$, and by general properties of summability, it is sufficient to prove that, for each $k\in \N$, the family $((\phi_n L_{\mathbf{e}}^n)(e_k))_{n\in \N}$ is summable in $V$.  Since $V$ is discrete and therefore complete,  it is sufficient to check that Cauchy's condition is satisfied. We may take $W:=\{0\}$ as a neighborhood  of zero in $V$. Let $J_{W}:=\{0,\cdots,k\}$. Because for every $n>k$, $L_{\mathbf{e}}^n(e_k)=0$, then $\displaystyle\sum_{n\in J}\left ( \phi_n L_{\mathbf{e}}^n)(e_k)\right )=0$ whenever $J$ is a finite subset of $I$ such that $J\cap J_W = \emptyset$. In what follows, the sum of a family $(\phi_n L^n_{\mathbf{e}})_{n\in\N}$ is the element of $\End(V)$ denoted by $\displaystyle\sum_{n\in \N}\phi_n L^n_{\mathbf{e}}$  where for every nonzero $v\in V$,
\begin{equation}
\displaystyle\left (\sum_{n\in \N}\phi_n L^n_{\mathbf{e}}\right )(v)=\sum_{n=0}^{\deg_{\mathbf{e}}(v)}\phi_n (L^n_{\mathbf{e}}(v))\ .
\end{equation}

We are now in a position to establish the main result concerning the expansion of any operator on $V$ in terms of ladder operators.
\begin{theoreme}[Endomorphism expansion in ladder operators]\label{mainthm}
Let $\mathbf{a}=(a_n)_{n\in \N}$ and $\mathbf{b}=(b_n)_{n\in \N}$ be two bases of   $V$ such that $b_0 \in \K a_0$; that is, there exists a nonzero scalar $\lambda:=\langle b_0,a_0\rangle$ such that $\lambda a_0 = b_0$. Then each $\phi\in\End(V)$ is the sum of the summable family $(P_n(R_{\mathbf{a}})L^n_{\mathbf{b}})_{n\in\N}$ where $(P_n)_{n\in\N}\in \K[\x]^{\N}$ is a sequence of polynomials that satisfies the following recursion equation
\begin{equation}
\left \{
\begin{array}{lll}
\lambda P_0(\mathbf{a})&=&\phi(b_0)\ ,\\
\lambda P_{n+1}(\mathbf{a})&=&\phi(b_{n+1})-\displaystyle\sum_{k=0}^n P_k(R_{\mathbf{a}})b_{n+1-k}\ .
\end{array}\right .
\end{equation}
(Note that due to Lemma~\ref{isomorphisme}  , for each $n\in\N$,   $P_n(\mathbf{a})$ uniquely defines $P_n\in \K[\x]$.)
\end{theoreme}

\begin{proof}
Since $\mathbf{b}$ is a basis, it is sufficient to prove that for each $n\in \N$,
\begin{equation}
\phi(b_n)=\left ( \displaystyle\sum_{k\in \N}P_k(R_{\mathbf{a}})L_{\mathbf{b}}^k\right )(b_n)\ .
\end{equation}
\begin{enumerate}
\item Case $n=0$:
\begin{equation}
\begin{array}{lll}
\left (\displaystyle \sum_{k\in \N}P_k(R_{\mathbf{a}})L_{\mathbf{b}}^k\right )(b_0)&=&P_0(R_{\mathbf{a}})(b_0)\\
&=&\lambda P_0(R_{\mathbf{a}})(a_0)\\
&=&\lambda P_0(\mathbf{a})\ (\mbox{according to Lemma~\ref{petitlem}})\\
&=&\phi(b_0)\ (\mbox{by assumption})\ .
\end{array}
\end{equation}
\item Case $n+1$, $n\in\N$:
\begin{equation}
\begin{array}{lll}
\left (\displaystyle \sum_{k\in \N}P_k(R_{\mathbf{a}})L_{\mathbf{b}}^k\right )(b_{n+1})&=&\displaystyle\sum_{k=0}^{n+1}P_k(R_{\mathbf{a}})b_{n+1-k}\\
&=&P_{n+1}(R_{\mathbf{a}})(b_0)+\displaystyle\sum_{k=0}^n P_k(R_{\mathbf{a}})b_{n+1-k}\\
&=&\lambda P_{n+1}(R_{\mathbf{a}})(a_0)+\displaystyle\sum_{k=0}^n P_k(R_{\mathbf{a}})b_{n+1-k}\\
&=&\lambda P_{n+1}(\mathbf{a})+\displaystyle\sum_{k=0}^n P_k(R_{\mathbf{a}})b_{n+1-k}\\
&=&\phi(b_{n+1})\ .
\end{array}
\end{equation}
\end{enumerate}
\end{proof}
%\begin{note}
%Without any change, we could replace the field $\K$ by an associative, commutative ring $R$ with a unit, and $V$ by a unital free $R$-module with a countable rank. Indeed every free modules with an infinite basis, even on a noncommutative ring, have a unique infinite rank (see for instance \cite{Dau94}, chap.~2, prop.~2-2.5).
%\end{note}

\begin{exemple}\label{evaluation_exmpl}
Suppose that $\K$ is a field of characteristic zero\footnote{The assumption on the characteristic of $\K$ is needed here because we consider denominators of the form $n!$.}. Consider $V:=\K[\x]$, $a_n := \x^n$ and $b_n := \frac{\x^n}{n!}$. Therefore $R_{\mathbf{a}}=X$, the operator of multiplication by $\x$; and, $L_{\mathbf{b}}=D$, the formal derivative of polynomials, which are the data of the classical result recalled in subsect.~\ref{sec3}. In Example~\ref{evaluation_exmpl}, we consider the functional $\epsilon \colon \K[\x]\rightarrow \K\subseteq \K[\x]$ that maps a polynomial to the sum of its coefficients. From Theorem~\ref{mainthm}, we know that $\epsilon=\displaystyle\sum_{n\geq 0}P_n(X)D^n$ and that 
\begin{equation}
P_{n+1}(\x)=\displaystyle \frac{1}{(n+1)!}-\displaystyle\sum_{k=0}^n P_k(\mathtt{x})\frac{\mathtt{x}^{n+1-k}}{(n+1-k)!}\ .
\end{equation} 
We can show by induction that $P_n(\x)=\frac{1}{n!}(1-\x)^n$, and then easily verify that $\epsilon=\displaystyle\sum_{n\geq 0}\frac{1}{n!}(1-X)^n D^n$ on the basis $\{\x^k\}_k$. Alternatively, we see that this operator is $\epsilon=e^{\mathtt{y} D}|_{\mathtt{y}=1-\x}\colon \x^n \mapsto (\x+\mathtt{y})^n|_{\mathtt{y}=1-\x}$. 
\end{exemple}

Let $\mathbf{a}=(a_n)_{n\in\N}$ and $\mathbf{b}=(b_n)_{n\in\N}$ be two bases of $V$.
Let us consider the following operators
\begin{equation}
L_{\mathbf{b},\boldsymbol{\beta}}b_n = \left \{
\begin{array}{lll}
0 & \mbox{if} & n=0\ ,\\
\beta_n b_{n-1} & \mbox{if} & n>0\ ,
\end{array}
\right .
\end{equation}
and
\begin{equation}\label{eq65}
R_{\mathbf{a},\boldsymbol{\alpha}}a_n = \alpha_n a_{n+1}
\end{equation}
where $\boldsymbol{\beta}:=(\beta_n)_{n\in\N}$, with $\beta_0:=1$,  and $\boldsymbol{\alpha}:=(\alpha_n)_{n\in\N}$ are sequences of nonzero scalars. These operators, which we may call respectively  {\emph{$\mathbf{b}$-relative lowering operator with coefficient sequence $\boldsymbol{\beta}$}} and {\emph{$\mathbf{a}$-relative raising operator with coefficient sequence $\boldsymbol{\alpha}$}}, seem to be straightforward  generalizations of the ladder operators as previously introduced; however, this is not entirely  the case. Actually, $L_{\mathbf{b},\boldsymbol{\beta}}$ and $R_{\mathbf{a},\boldsymbol{\alpha}}$ are respectively equal to
some ``usual" ladder operators $L_{\boldsymbol{\beta}^{-1}\cdot\mathbf{b}}$ and $R_{\boldsymbol{\alpha}\cdot\mathbf{a}}$ where $\boldsymbol{\beta}^{-1}\cdot\mathbf{b}:=(b_n')_{n\in\N}$ with $b_n'=\left(\displaystyle\prod_{i=0}^n\beta_i\right)^{-1} b_n$ (resp. $\boldsymbol{\alpha}\cdot\mathbf{a}=(a_n')_{n\in\N}$ where $a_n'=\displaystyle(\prod_{i=0}^{n-1}\alpha_i)a_n$ for $n>0$, and $a_0'=a_0$). If $b_0' \in \K a_0'$ (or, equivalently, if $b_0\in\K a_0$, because $b_0'=\frac{b_0}{\beta_0}=b_0$ and $a_0'=a_0$), then we can apply Theorem~\ref{mainthm} with the operators $L_{\mathbf{b},\boldsymbol{\beta}}$ and $R_{\mathbf{a},\boldsymbol{\alpha}}$, just by replacing $\mathbf{a}$ by
$\boldsymbol{\alpha}\cdot\mathbf{a}$, $\mathbf{b}$ by $\boldsymbol{\beta}^{-1}\cdot\mathbf{b}$. When $\mathbf{a}=\mathbf{b}$, we say that $L_{\mathbf{a},\boldsymbol{\beta}}$ and $R_{\mathbf{a},\boldsymbol{\alpha}}$ are {\emph{$\mathbf{a}$-relative ladder operators}}  with coefficients $\boldsymbol{\beta}$ and $\boldsymbol{\alpha}$ respectively. Such a pair of operators - used in the following subsection -  satisfy the  rather general commutation rule
\begin{equation}
D_{\mathbf{a},\boldsymbol{\beta},\boldsymbol{\alpha}}:=[L_{\mathbf{a},\boldsymbol{\beta}},R_{\mathbf{a},\boldsymbol{\alpha}}]=L_{\mathbf{a},\boldsymbol{\beta}}R_{\mathbf{a},\boldsymbol{\alpha}} - R_{\mathbf{a},\boldsymbol{\alpha}}L_{\mathbf{a},\boldsymbol{\beta}}
\end{equation}
where $D_{\mathbf{a},\boldsymbol{\beta},\boldsymbol{\alpha}}$ is the operator defined by
\begin{equation}
D_{\mathbf{a},\boldsymbol{\beta},\boldsymbol{\alpha}} a_n = \left \{
\begin{array}{lll}
(\alpha_0 \beta_1)a_0 & \mbox{if} & n=0\ ,\\
(\alpha_n \beta_{n+1}-\alpha_{n-1}\beta_n)a_{n} & \mbox{if} & n>0\ ,
\end{array}
\right .
\end{equation}
which we call the {\emph{diagonal operator}} associated with $L_{\mathbf{a},\boldsymbol{\beta}}$ and $R_{\mathbf{a},\boldsymbol{\alpha}}$.
\begin{note}
It is possible to define a similar $D_{\mathbf{b},\mathbf{a}}\in\End(V)$ associated with any ladder operators $L_{\mathbf{b}}$ and $R_{\mathbf{a}}$ by $D_{\mathbf{b},\mathbf{a}} := [L_{\mathbf{b}},R_{\mathbf{a}}]$, which defines the commutation relation between $L_{\mathbf{b}}$ and $R_{\mathbf{a}}$. (In particular, $D_{\mathbf{a},\boldsymbol{\beta},\boldsymbol{\alpha}} = D_{\boldsymbol{\beta}^{-1}\cdot\mathbf{a},\boldsymbol{\alpha}\cdot\mathbf{a}}$.) Furthermore, when the two bases $\mathbf{a}$ and $\mathbf{b}$ are related by $b_0\in\K a_0$ as in Theorem~\ref{mainthm}, then, as an operator on $V$, $D_{\mathbf{b},\mathbf{a}}$ is the sum of a summable family $(P_n(R_{\mathbf{a}})L_{\mathbf{b}}^n)_{n\in\N}$, and therefore the commutation relation is given by
\begin{equation}
L_{\mathbf{b}}R_{\mathbf{a}}=R_{\mathbf{a}}L_{\mathbf{b}}+\displaystyle\sum_{n\in\N}P_n(R_{\mathbf{a}})L^{n}_{\mathbf{b}}\ .
\end{equation}
\end{note}

\subsection{Extension to formal infinite linear combinations}

\subsubsection{Preliminaries: topology and duality}

Let $\K$ be a field (of any characteristic). Let $V$ be a countable-dimensional $\K$-vector space, and $\mathbf{e}:=(e_n)_{n\in\N}$ be a basis of $V$. The vector space $V$ can be considered  as the $\N$-graded vector space $V_{\mathbf{e}}:=\displaystyle\bigoplus_{n\in\N}\K e_n$. There exists a natural decreasing filtration associated with this grading which is defined by $V=V_{\mathbf{e}}=\displaystyle\bigcup_{n\in\N}F_n(V_{\mathbf{e}})$ where $F_n(V_{\mathbf{e}}):=\displaystyle\bigoplus_{k\geq n}\K e_k$. This filtration is separated, \emph{i.e.}, $\displaystyle\bigcap_{n\in\N}F_n(V_{\mathbf{e}})=(0)$. Now suppose that $\K$ has the discrete topology. The subsets $F_n(V_{\mathbf{e}})$ define a fundamental system of neighbourhoods of zero of a Hausdorff $\K$-vector topology on $V=V_{\mathbf{e}}$ (see~\cite{Bou-CA}). This (metrizable) topology may be equivalently described in terms of an order function. Define $\omega_{\mathbf{e}} : V_{\mathbf{e}} \rightarrow \N\cup\{+\infty\}$ by
\begin{equation}
\omega_{\mathbf{e}}(v)=\left \{
\begin{array}{lll}
\min\{n\in\N : \langle v,e_n\rangle\not=0\}&\mbox{if} & v\not=0\ ,\\
+\infty & \mbox{if} & v=0
\end{array}
\right .
\end{equation}
for $v\in V$.
The completion $\widehat{V}_{\mathbf{e}}$ of $V_{\mathbf{e}}$ for this topology is canonically identified with the $\K$-vector space $\displaystyle\prod_{n\in\N}\K e_n$ - that is, the set of all families $(v_n)_{n\in\N}$ with $v_n\in\K e_n$ for each integer $n$ -  equipped with the product topology of discrete topologies on each factor $\K e_n$. Each element $S$ of $\widehat{V}_{\mathbf{e}}$ may be uniquely seen as a formal infinite linear combination $S=\displaystyle\sum_{n\in \N}\langle S,e_n\rangle e_n$, where $\langle S,e_n\rangle e_n = v_n$ and $S=(v_n)_{n\in\N}$ (it is not difficult to prove that the family $(\langle S,e_n\rangle e_n)_{n\in\N}$ is actually summable). The topology induced by $\widehat{V}_{\mathbf{e}}$ on $V_{\mathbf{e}}$ is the same as the topology defined by the filtration.  The order function is extended to $\widehat{V}_{\mathbf{e}}$ by
\begin{equation}
\omega_{\mathbf{e}}(S)=\left \{
\begin{array}{lll}
\min\{n\in\N : \langle S,e_n\rangle\not=0\} & \mbox{if} & S\not=0\ ,\\
+\infty & \mbox{if} & S=0
\end{array}
\right .
\end{equation}
for $S\in\widehat{V}_{\mathbf{e}}$, and may be used to describe the topology of the completion. For instance, a sequence $(S_n)_{n\in\N}$ of formal infinite linear combinations converges to zero if, and only if, $\displaystyle\lim_{n\rightarrow\infty}\omega_{\mathbf{e}}(S_n)=+\infty$;in other terms, for every $n\in \N$ there are only finitely many $k\in \N$ such that $\langle S_k, e_n\rangle\not=0$. This topology is sometimes referred to s the {\emph{formal topology}} (see~\cite{Chen, Kac}), and, $\widehat{V}_{\mathbf{e}}$ is then the {\emph{formal completion}} of the $\N$-graded vector space $V_{\mathbf{e}}:=\displaystyle\bigoplus_{n\in\N}\K e_n$.

\begin{note}
If $\mathbf{a}:=(a_n)_{n\in\N}$ and $\mathbf{b}:=(b_n)_{n\in\N}$ are two bases of $V$, then the isomorphism $\Phi$ of $V$ that maps $a_n$ to $b_n$ for each $n\in \N$  is also a homeomorphism from $V_{\mathbf{a}}$ to $V_{\mathbf{b}}$ considered  as spaces equipped with their respective filtrations. It turns out that $\Phi$ may be extended to a homeomorphism $\widehat{\Phi}$ from $\widehat{V}_{\mathbf{a}}$ to $\widehat{V}_{\mathbf{b}}$. Although the two spaces are homeomorphic, we cannot canonically identify them. Indeed, let us consider the sequence $\mathbf{b}:=(b_n)_{n\in\N}$ defined by $b_n := \displaystyle\sum_{k=0}^n a_k$, where $\mathbf{a}:=(a_n)_{n\in\N}$ is another basis. Then $\mathbf{b}$ is a basis of $V$: suppose that for some $n\in\N$, we have $\displaystyle\sum_{i=0}^n \alpha_i b_i=0$ with $\alpha_i\in\K$. Then $\displaystyle\sum_{i=0}^{n}\alpha_i\left ( \sum_{k=0}^i a_k\right )=0$ which is equivalent to $\displaystyle(\sum_{i=0}^n\alpha_i)b_0+(\sum_{i=1}^n\alpha_i)b_1+\cdots+(\alpha_{n-1}+\alpha_n)a_{n-1} + \alpha_n a_n=0$. Then $\alpha_i=0$ for every $i=0,\cdots,n$, and $\{b_i : i=0,\cdots,n\}$ is linearly independent. Using the classical M\"obius inversion, we obtain
\begin{equation}
a_n = \left \{
\begin{array}{lll}
b_0 & \mbox{if} & n=0\ ,\\
b_n - b_{n-1} & \mbox{if} & n>0
\end{array}
\right .
\end{equation}
which proves that $V$ is generated by $\mathbf{b}$.
Now we have  $\displaystyle\lim_{n\rightarrow \infty}b_n=0$ in the topology of $V_{\mathbf{b}}$  , but $\displaystyle\lim_{n\rightarrow \infty}b_n = \sum_{n=0}^{\infty}a_n$ in $\widehat{V}_{\mathbf{a}}$. (Note however that
$\displaystyle\lim_{n\rightarrow\infty}a_n=0$ also in $V_{\mathbf{b}}$ because $\omega_{\mathbf{b}}(a_{n})=n-1$ for every $n\in \N\setminus\{0\}$, and then $\displaystyle\lim_{n\rightarrow\infty}\omega_{\mathbf{b}}(a_n)=+\infty$.) The problem is due to the fact that the order function depends on the choice of the basis. %Now if the two bases $\mathbf{a}=(a_n)_{n\in\N}$ and $\mathbf{b}=(b_n)_{n\in\N}$ are linked by $a_n = \langle a_n,b_n\rangle b_n$, with $\langle a_n,b_n\rangle \not=0$, for every $n\in \N$, {\it i.e.}, $a_n \in\K b_n$ for each $n\in\N$, then $\omega_{\mathbf{a}}(v)=\omega_{\mathbf{b}}(v)$ for every
%$v\in V_{\mathbf{a}}=V$. (Note also that $\deg_{\mathbf{a}}(v)=\deg_{\mathbf{b}}(v)$ for every $v\in V_{\mathbf{a}}=V$.) In other terms, $\Id_V$ is an homeomorphism from $V$ with the topology associated with $V_{\mathbf{a}}$ to $V$ with the topology associated with $V_{\mathbf{b}}$, and both topologies are equivalent. In this case, $\widehat{V}_{\mathbf{a}}$ and $\widehat{V}_{\mathbf{b}}$ are ``essentially" the same.
\end{note}

 We now introduce the {\emph{(duality) pairing}} $\langle .|.\rangle : V_{\mathbf{e}}\times \widehat{V}_{\mathbf{e}} \rightarrow \K$ defined by $\langle P|S\rangle := \displaystyle\sum_{n=0}^{\deg_{\mathbf{e}}(P)}\langle P,e_n\rangle \langle S,e_n\rangle$, for $P\in V_{\mathbf{e}}$ and $S\in\widehat{V}_{\mathbf{e}}$. This pairing, also considered in \cite{Stanley2}, satisfies in particular
\begin{equation}
\langle e_i|e_j\rangle = \langle e_i,e_j\rangle = \delta_{i,j}= \left \{
\begin{array}{lll}
0 & \mbox{if} & i\not=j\ ,\\
1 & \mbox{if} & i=j
\end{array}
\right .
\end{equation}
for each $i,j\in\N$, and more generally $\langle P|e_j\rangle=\langle P,e_j\rangle$, $\langle e_i|S\rangle=\langle S,e_i\rangle$ for every $P\in V_{\mathbf{e}}$, $S\in\widehat{V}_{\mathbf{e}}$.

The algebraic dual space $V_{\mathbf{e}}^*$ of $V_{\mathbf{e}}$ is isomorphic to $\widehat{V}_{\mathbf{e}}$. Indeed let $\ell\in V_{\mathbf{e}}^*$ and define $S_{\ell}:=\displaystyle\sum_{n\in\N}\ell(e_n)e_n\in \widehat{V}_{\mathbf{e}}$. Then
$\ell(P)=\langle P|S_{\ell}\rangle$. The linear mapping $\ell\mapsto S_{\ell}$ is clearly one-to-one. It is also onto because for each $S\in \widehat{V}_{\mathbf{e}}$, $P\mapsto \langle P,S\rangle$ is easily seen as a linear form over $V_{\mathbf{e}}$.

The topological dual space $\widehat{V}_{\mathbf{e}}'$ of $\widehat{V}_{\mathbf{e}}$ is isomorphic to $V_{\mathbf{e}}$. Indeed let us consider a linear continuous form $\ell$ of $\widehat{V}_{\mathbf{e}}$. Since $\ell$ is continous, for every $S\in \widehat{V}_{\mathbf{e}}$, $\ell(S)=\displaystyle \sum_{n\geq 0}\langle S,e_n\rangle \ell(e_n)$ and the sum is convergent in $\K$ discrete. Therefore there is an integer $N$ such that for every $n\geq N$, $\langle S,e_n\rangle\ell(e_n)=0$. If we choose $S=\displaystyle \sum_{n\geq 0}e_n$, then it means that for $n$ large enough, $\ell(e_n)=0$. Then $P_{\ell}=\displaystyle\sum_{n\geq 0}\ell(e_n)e_n$ is actually an element of $V_{\mathbf{e}}$ which satisfies $\langle P_{\ell}|S\rangle = \ell(S)$ for every formal infinite linear combination $S$.
Now suppose that $P_{\ell}=0$ for $\ell\in\widehat{V}_{\mathbf{e}}'$. Then for every $n\in\N$, $\ell(e_n)=\langle P_{\ell}|e_n\rangle=\langle P_{\ell},e_n\rangle=0$. The linear form is null on the dense subset $V_{\mathbf{e}}$ of $\widehat{V}_{\mathbf{e}}$, and, by continuity, $\ell$ is also equal to zero on the closure. Let $P\in V_{\mathbf{e}}$. Then $\ell:=S\mapsto \langle P|S\rangle$ is a linear form on $\widehat{V}_{\mathbf{e}}$ such that $P_{\ell}=P$. Moreover,
$\ell$ is clearly continuous. In summary, the pairing performs the following isomorphisms.
\begin{equation}
\begin{array}{lll}
V_{\mathbf{e}}^* &\cong& \widehat{V}_{\mathbf{e}}\ ,\\
\widehat{V}_{\mathbf{e}}'&\cong& V_{\mathbf{e}}\ .
\end{array}
\end{equation}
The respective isomorphisms are given by
\begin{equation}
\Phi : V_{\mathbf{e}}^{*} \rightarrow \widehat{V}_{\mathbf{e}}
\end{equation}
and
\begin{equation}
\Psi : \widehat{V}_{\mathbf{e}}' \rightarrow V_{\mathbf{e}}
\end{equation}
such that for every $P\in V_{\mathbf{e}}$, $S\in \widehat{V}_{\mathbf{e}}$, if $\ell\in V_{\mathbf{e}}^{*}$, then
\begin{equation}
\langle P|\Phi(\ell)\rangle = \ell(P)
\end{equation}
while
\begin{equation}
\Phi^{-1}(S)(P)=\langle P|S\rangle
\end{equation}
and if $\ell\in \widehat{V}_{\mathbf{e}}'$, then
\begin{equation}
\langle \Psi(\ell)|S\rangle = \ell(S)
\end{equation}
and
\begin{equation}
\Psi^{-1}(P)(S)=\langle P|S\rangle\ .
\end{equation}
We may use these isomorphisms to define the natural notion of {\em transpose} in this setting. The {\em{transpose}} of $\phi\in\End(V_{\mathbf{e}})$ is $\phi^{t}\in\End(\widehat{V}_{\mathbf{e}}^*)$ such that for every $P\in V_{\mathbf{e}}$ and every $S\in \widehat{V}_{\mathbf{e}}$, $\langle \phi P|S\rangle = \langle P|\phi^{t}S\rangle$. Actually, $\phi^{t}$ is defined as
\begin{equation}
\begin{array}{llll}
\phi^{t}:& \widehat{V}_{\mathbf{e}} & \rightarrow & \widehat{V}_{\mathbf{e}}\\
& S &\mapsto & \Phi(\Phi^{-1}(S)\circ\phi)\ .
\end{array}
\end{equation}
Indeed, for every $P\in V_{\mathbf{e}}$, the following holds.
\begin{equation}
\begin{array}{lll}
\langle P|\phi^{t}(S)\rangle&=&\langle P|\Phi(\Phi^{-1}(S)\circ\phi\rangle\\
&=&(\Phi^{-1}(S))(\phi(P))\\
&=&\langle\phi(P)|S\rangle\ .
\end{array}
\end{equation}
By duality, it is also possible to define a transpose for $\phi\in\End(\widehat{V}_{\mathbf{e}})$ but continuity has to be taken into account. Indeed,
let $\phi\in\End(\widehat{V}_{\mathbf{e}})$ be a continuous endomorphism. We can define ${}^t\phi\in\End(V_{\mathbf{e}})$ by
\begin{equation}
{}^t\phi(P):=\Psi(\Psi^{-1}(P)\circ\phi)
\end{equation}
for every $P\in V_{\mathbf{e}}$. Note that since $\phi$ is continuous (and linear), $\Psi^{-1}(P)\circ\phi\in \widehat{V}_{\mathbf{e}}'$. Then, for every $P\in V_{\mathbf{e}}$ and $S\in \widehat{V}_{\mathbf{e}}$, we have
\begin{equation}
\langle P|\phi(S)\rangle = \langle {}^t\phi(P)|S\rangle\ .
\end{equation}
Indeed,
\begin{equation}
\begin{array}{lll}
\langle {}^t\phi(P)|S\rangle&=&\langle\Psi(\Psi^{-1}(P)\circ\phi)|S\rangle\\
&=&(\Psi^{-1}(P))(\phi(S))\\
&=&\langle P|\phi(S)\rangle\ .
\end{array}
\end{equation}
\begin{lemme}\label{transppolyendoiscontinuous}
For each $\phi\in\End(V_{\mathbf{e}})$, $\phi^t$ is a continous endomorphism of $\widehat{V}_{\mathbf{e}}$. Moreover, $\phi={}^t(\phi^{t})$. Dually, for every continuous endomorphism $\phi$ of $\widehat{V}_{\mathbf{e}}$, $\phi=({}^t\phi)^t$.
\end{lemme}
\begin{proof}
Let $\phi\in \End(V_{\mathbf{e}})$ and $\{S_n\}_n$ be a sequence of infinite linear combinations that converges to zero. Let $k\in \N$. By definition ot the transpose, $\langle\phi^t(S_n),e_k\rangle=\displaystyle\sum_{i\geq 0}\langle\phi(e_k),e_i\rangle\langle S_n,e_i\rangle$. Since $S_n\rightarrow 0$, for every $i$, there is $N_i$ such that for all $n\geq N_i$, $\langle S_n,e_i\rangle=0$. Therefore we can find $N_k$ such that $n\geq N_k$ implies $\langle S_n,e_i\rangle=0$ for every $i\leq \deg_{\mathbf{e}}(\phi(e_k))$, and then for such $n$, $\langle \phi^t(S_n),e_k\rangle=0$, so $\phi^t(S_n)\rightarrow 0$, and $\phi^t$ is continuous. Now let us prove that $\phi={}^t(\phi^t)$. For every $P,S$, we have $\langle \phi(P)|S\rangle=\langle P|\phi^t(S)\rangle=\langle {}^t(\phi^t)(P)|S\rangle$ (the second equality is valid since $\phi^{t}$ is continuous). Therefore for every $i,j$, $\langle\phi(e_i),e_j\rangle=\langle \phi(e_i)|e_j\rangle=\langle {}^t(\phi^t)(e_i)|e_j\rangle=\langle{}^t(\phi^t)(e_i),e_j\rangle$ which is sufficient to prove the expected equality. Finally, let $\phi$ be a continuous endomorphism of $\widehat{V}_{\mathbf{e}}$. For every $P,S$, one has $\langle P|\phi(S)\rangle=\langle {}^t\phi(P)|S\rangle=\langle P|({}^t \phi)^t(S)\rangle$, and in particular for every $i$, 
$\langle \phi(S),e_i\rangle=\langle e_i|\phi(S)\rangle=\langle e_i|({}^t \phi)^t(S)\rangle=\langle({}^t \phi)^t(S),e_i\rangle$, which proves that $\phi(S)=({}^t\phi)^t(S)$ (by definition of $\widehat{V}_{\mathbf{e}}$). 
\end{proof}

Let $\mathbf{a}$ and $\mathbf{b}$ be two bases of $V$. Let $L_{\mathbf{b},\boldsymbol{\beta}}$ (resp. $R_{\mathbf{a},\boldsymbol{\alpha}}$) be a $\mathbf{b}$-relative lowering operator (resp. $\mathbf{a}$-relative raising operator) with coefficient sequence $\mathbf{\beta}=(\beta_n)_{n\in\N}$ with $\beta_0=1$ (resp. $\mathbf{\alpha}=(\alpha_n)_{n\in\N}$). These operators are clearly continuous on $V_{\mathbf{b}}$ (resp. on $V_{\mathbf{a}}$), and therefore extend uniquely as continuous endomorphisms of
the completions $\widehat{V}_{\mathbf{b}}$ and $\widehat{V}_{\mathbf{a}}$. Their respective extensions $\widehat{L}_{\mathbf{b},\boldsymbol{\beta}}$ and $\widehat{R}_{\mathbf{a},\boldsymbol{\alpha}}$ are precisely defined by
\begin{equation}\label{extensionofL}
\widehat{L}_{\mathbf{b},\boldsymbol{\beta}}(S)=\displaystyle\sum_{n\geq 0}\langle S,b_n\rangle L_{\mathbf{b},\boldsymbol{\beta}}b_n=\sum_{n\geq 1}\langle S,b_n\rangle \beta_n b_{n-1}
=\sum_{n\geq 0}\langle S,b_{n+1}\rangle \beta_{n+1} b_n
\end{equation}
and
\begin{equation}
\widehat{R}_{\mathbf{a},\boldsymbol{\alpha}}(S)=\displaystyle \sum_{n\geq 0}\langle S,a_n\rangle R_{\mathbf{a},\boldsymbol{\alpha}}a_n=\sum_{n\geq 0}\langle S,a_n\rangle \alpha_n a_{n+1}=
\sum_{n\geq 1}\langle S,a_{n-1}\rangle\alpha_{n-1}a_n\ .
\end{equation}
They correspond to the operators $D$ and $U$ of \cite{Stanley2} associated with the graded (locally finite) posets
$b_0 \rightarrow b_1 \rightarrow \frac{1}{\beta_1}b_2 \rightarrow \frac{1}{\beta_1\beta_2}b_3\rightarrow\cdots$ and
$a_0 \rightarrow \alpha_0 a_1 \rightarrow \alpha_0\alpha_1 a_2\rightarrow \alpha_0\alpha_1\alpha_2 a_3 \rightarrow\cdots$

We may use the duality pairing in order to find the transpose mappings of both $L_{\mathbf{b},\boldsymbol{\beta}}$ and $R_{\mathbf{a},\boldsymbol{\alpha}}$.
\begin{lemme}\label{transpose-raising}
Let $R_{\mathbf{a},\boldsymbol{\alpha}}$ be the $\mathbf{a}$-relative raising operator with coefficient sequence $\boldsymbol{\alpha}=(\alpha_n)_{n\in\mathbb{N}}$.
The transpose of $R_{\mathbf{a},\boldsymbol{\alpha}}$ is the extension $\widehat{L}_{\mathbf{a},\boldsymbol{\gamma}}$ to the completion $\widehat{V}_{\mathbf{a}}$ of the $\mathbf{a}$-relative lowering operator $L_{\mathbf{a},\boldsymbol{\alpha}\downarrow}$ with coefficient sequence $\boldsymbol{\alpha}\downarrow:=(\gamma_n)_{n\in\N}$ where
\begin{equation}
\gamma_n := \left \{
\begin{array}{lll}
1 & \mbox{if} & n=0\ ,\\
\alpha_{n-1} & \mbox{if} & n> 0\ .
\end{array}
\right .
\end{equation}
\end{lemme}

\begin{proof}
Let $n\in\N$ and $S\in \widehat{V}_{\mathbf{a}}$. According to Equation~(\ref{eq65}), $\langle R_{\mathbf{a},\boldsymbol{\alpha}}a_n|S\rangle=\alpha_n\langle a_{n+1}|S\rangle=\alpha_n\langle S,a_{n+1}\rangle=\displaystyle\left \langle a_n|\sum_{k\geq 0}\langle S,a_{k+1}\rangle \alpha_k a_k\right \rangle=\langle a_n|\widehat{L}_{\mathbf{a},\boldsymbol{\alpha}\downarrow}\rangle$ (the last equality comes from Equation~\ref{extensionofL}). Multiplying both (leftmost and rightmost) sides  with $\langle P|a_n\rangle$ (for some $P\in V_{\mathbf{a}}$) and summing over
$n$ gives the result.
\end{proof}

\begin{lemme}\label{transpose-lowering}
Let $L_{\mathbf{b},\boldsymbol{\beta}}$ be the $\mathbf{b}$-relative lowering operator with coefficient sequence $\boldsymbol{\beta}=(\beta_n)_{n\in\N}$. The transpose $L_{\mathbf{b},\boldsymbol{\beta}}^{t}$ of $L_{\mathbf{b},\boldsymbol{\beta}}$ is the extension $\widehat{R}_{\mathbf{b},\boldsymbol{\beta}\uparrow}$ to $\widehat{V}_{\mathbf{b}}$ of the $\mathbf{b}$-relative raising operator $R_{\mathbf{b},\boldsymbol{\beta}\uparrow}$ with coefficient sequence $\boldsymbol{\beta}\uparrow:=(\gamma_n)_{n\in\N}$, where for each $n\in\N$, $\gamma_n := \beta_{n+1}$.
\end{lemme}
\begin{proof}
This proof is so similar to the proof of Lemma~\ref{transpose-raising}, that it can be
omitted.
\end{proof}

It is also possible to determine the transpose of the extension of the ladder operators to the completion $\widehat{V}_{\mathbf{e}}$. Several lemmas are given below to answer this question. The first one does not need a proof. 
\begin{lemme}\label{updown-and-downup}
Let $\boldsymbol{\beta}=(\beta_n)_{n\in\N}$ be any sequence of elements of $\K$ such that $\beta_0=1$. We have
\begin{equation}
\boldsymbol{\beta}=\boldsymbol{\beta}\uparrow\downarrow\ .
\end{equation}
Let $\boldsymbol{\alpha}=(\alpha_n)_{n\in\N}$ be any sequence of elements of $\K$. We have
\begin{equation}
\boldsymbol{\alpha}=\boldsymbol{\alpha}\downarrow\uparrow\ .
\end{equation}
\end{lemme}
%\begin{proof}
%Let $\boldsymbol{\beta}\uparrow:=(\gamma_n)_{n\in \N}=\boldsymbol{\gamma}$ with $\gamma_n = \beta_{n+1}$ for each integer $n$. Let $(\delta_n)_{n\in\N}=\boldsymbol{\gamma}\downarrow=\boldsymbol{\beta}\uparrow\downarrow$.  By definition, $\delta_0 = 1$ and for each $n\in\N$, $\delta_n = \gamma_{n-1}=\beta_n$. Then for every $n\in\N$, $\beta_n = \delta_n$.
%\end{proof}
\begin{lemme}\label{tranpose-op-series}
Let $\mathbf{e}=(e_n)_{n\in\N}$ be a basis of $V$. Let $\boldsymbol{\beta}=(\beta_n)_{n\in\N}$ be a sequence of nonzero scalars such that $\beta_0=1$, and $\boldsymbol{\alpha}=(\alpha_n)_{n\in\N}$ be any sequence of nonzero scalars. Then we have
\begin{equation}
{}^t\widehat{L}_{\mathbf{e},\boldsymbol{\beta}}=R_{\mathbf{e},\boldsymbol{\beta}\uparrow}\ \mbox{and}\ {}^t\widehat{R}_{\mathbf{e},\boldsymbol{\alpha}}=L_{\mathbf{e},\boldsymbol{\alpha}\downarrow}\ .
\end{equation}
\end{lemme}
\begin{proof}
Since $\widehat{L}_{\mathbf{e},\boldsymbol{\beta}}$ and $\widehat{R}_{\mathbf{e},\boldsymbol{\alpha}}$ are continuous endomorphisms of $\widehat{V}_{\mathbf{e}}$ they admit transposes which are endomorphisms of $V_{\mathbf{e}}$. According to lemmas~\ref{transpose-raising} and~\ref{updown-and-downup}, $R^t_{\mathbf{e},\boldsymbol{\beta}\uparrow}=\widehat{L}_{\mathbf{e},\boldsymbol{\beta}\uparrow\downarrow}=\widehat{L}_{\mathbf{e},\boldsymbol{\beta}}$. Then, ${}^t\widehat{L}_{\mathbf{e},\boldsymbol{\beta}}={}^t(R^t_{\mathbf{e},\boldsymbol{\beta}\uparrow})=R_{\mathbf{e},\boldsymbol{\beta}\uparrow}$ (according to lemma~\ref{transppolyendoiscontinuous}). The case of ${}^t\widehat{R}_{\mathbf{e},\boldsymbol{\alpha}}$ is treated in a similar way.
\end{proof}
%\begin{proof}
%Let $\mathbf{\gamma}=(\gamma_n)_{n\in\N}=\boldsymbol{\alpha}\downarrow$, \emph{i.e.}, $\gamma_0=1$ and $\gamma_{n+1}=\alpha_n$. Let also define
%$(\delta_n)_{n\in\N}=\boldsymbol{\gamma}\uparrow=\boldsymbol{\alpha}\downarrow\uparrow$. Then, $\delta_n=\gamma_{n+1}=\alpha_n$ for every $n\in\N$, and so $\boldsymbol{\alpha}=\boldsymbol{\alpha}\downarrow\uparrow$.
%\end{proof}

\subsubsection{Extension of Theorem~\ref{mainthm} to formal infinite linear combinations}

In what follows, our intention is to generalize Theorem~\ref{mainthm} to the case of continuous endomorphisms on formal infinite linear combinations. To this end, we suppose that $\widehat{V}_{\mathbf{e}}$ is equipped with the $V_{\mathbf{e}}$-weak topology, that is, the weakest topology for which the mappings $\Psi^{-1}(P):S\in \widehat{V}_{\mathbf{e}} \mapsto \langle P|S\rangle\in \K$, defined for a given $P\in V_{\mathbf{e}}$, are continuous. Since $V_{\mathbf{e}}$ is isomorphic to $\widehat{V}_{\mathbf{e}}'$ (when $\widehat{V}_{\mathbf{e}}$ is equipped with its formal topology previously introduced), it is the so-called weak-$*$ topology. This topology turns $\widehat{V}_{\mathbf{e}}$ into a Hausdorff topological space (with $\K$ discrete). It is obvious that the duality pairing $\langle.|.\rangle$ is separately continuous on $V_{\mathbf{e}}\times \widehat{V}_{\mathbf{e}}$ where $V_{\mathbf{e}}$ is discrete and $\widehat{V}_{\mathbf{e}}$ has the $V_{\mathbf{e}}$-weak topology. Thus, a family $(S_i)_{i\in I}\in \widehat{V}_{\mathbf{e}}^I$ is summable whenever for every $P\in V_{\mathbf{e}}$, the family $(\langle P|S_i\rangle)_{i\in I}$ is summable in $\K$, and, in this case, $\langle P|\displaystyle\sum_{i\in I}S_i\rangle =
\sum_{i\in I}\langle P|S_i\rangle$.\\

Now suppose that the vector space of continuous endomorphisms of $\widehat{V}_{\mathbf{e}}$ has the topology of simple convergence. (We also suppose the same for $\End(V_{\mathbf{e}})$, with $V_{\mathbf{e}}$ equipped with the discrete topology.) In this particular topology, each family of continuous endomorphisms $(\widehat{R}_{\mathbf{e},\boldsymbol{\alpha}}^n\phi_n)_{n\in\N}$ in $\End(\widehat{V}_{\mathbf{e}})^{\N}$, where $\phi_n$ is a continuous endomorphism of $\widehat{V}_{\mathbf{e}}$ for each integer $n$, is a summable family. In order to check this, let
$P\in V_{\mathbf{e}}$ and $S\in \widehat{V}_{\mathbf{e}}$. We have ${}^t(\widehat{R}_{\mathbf{e},\boldsymbol{\alpha}}^n \phi_n)={}^t \phi_n L^{n}_{\mathbf{e},\boldsymbol{\alpha}\downarrow}\in \End(V_{\mathbf{e}})$. The family $({}^t \phi_n L^{n}_{\mathbf{e},\boldsymbol{\alpha}\downarrow})_{n\in\N}$ is summable in $\End(V_{\mathbf{e}})$, and we have
\begin{equation}
\begin{array}{lll}
\langle \displaystyle\sum_{n\in\N}{}^t \phi_n L^{n}_{\mathbf{e},\boldsymbol{\alpha}\downarrow}(P)|S\rangle&=&
\displaystyle\langle \sum_{n=0}^{\deg_{\mathbf{e}}(P)}{}^t \phi_n L^{n}_{\mathbf{e},\boldsymbol{\alpha}\downarrow}(P)|S\rangle\\
&=&\displaystyle\sum_{n=0}^{\deg_{\mathbf{e}}(P)}\langle {}^t \phi_n L^{n}_{\mathbf{e},\boldsymbol{\alpha}\downarrow}(P)|S\rangle\\
&=&\displaystyle \sum_{n=0}^{\deg_{\mathbf{e}}(P)}\langle P|\widehat{R}_{\mathbf{e},\boldsymbol{\alpha}}^n \phi_n S\rangle\\
&=&\displaystyle\langle P|\sum_{n=0}^{\deg_{\mathbf{e}}(P)}\widehat{R}_{\mathbf{e},\boldsymbol{\alpha}}^n \phi_n S\rangle\ .
\end{array}
\end{equation}
Moreover for every $m>\deg_{\mathbf{e}}(P)$,
\begin{equation}
\langle P|\displaystyle\sum_{n=\deg_{\mathbf{e}}(P)}^{m}\widehat{R}_{\mathbf{e},\boldsymbol{\alpha}}^n \phi_n S\rangle = \langle \sum_{n=\deg_{\mathbf{e}}(P)}^{m}{}^t \phi_n L^{n}_{\mathbf{e},\boldsymbol{\alpha}\downarrow}(P)|S\rangle=0\ .
\end{equation}
Therefore, we obtain a summable series in $\K$ discrete, and so is $(\widehat{R}_{\mathbf{e},\boldsymbol{\alpha}}^n\phi_n)_{n\in\N}$.

The generalization of Theorem~\ref{mainthm} to the case of continuous operators on formal infinite linear combinations is given below.
\begin{theoreme}\label{series-mainthm}
Let $\mathbf{\alpha}=(\alpha_n)_{n\in\N}$ be any sequence of nonzero scalars, and $\mathbf{\beta}=(\beta_n)_{n\in\N}$ be a sequence of nonzero scalars with $\beta_0=1$. Let $\phi$ be any continuous element of $\End(\widehat{V}_{\mathbf{e}})$. Then there exists a sequence of polynomials $(P_n)_{n\in\N}\in\K[\x]^{\mathbb{N}}$ such that $\phi$ is equal to the sum of the summable family $(\widehat{R}_{\mathbf{e},\boldsymbol{\beta}\uparrow}^n P_n(\widehat{L}_{\mathbf{e},\boldsymbol{\alpha}\downarrow}))_{n\in\N}$.
\end{theoreme}
\begin{proof}
By Theorem~\ref{mainthm}, ${}^t\phi = \displaystyle\sum_{n\in\N}P_n(R_{\mathbf{e},\boldsymbol{\alpha}})L_{\mathbf{e},\boldsymbol{\beta}}^n$ (sum of a summable family). Then, using the duality pairing, we check that $\phi = \displaystyle\sum_{n\in\N}\widehat{R}_{\mathbf{e},\boldsymbol{\beta}\uparrow}^nP_n(\widehat{L}_{\mathbf{e},\boldsymbol{\alpha}\downarrow})$ (sum of a summable family).
\end{proof}
\begin{corollaire}
Under the same assumptions as those of Theorem~\ref{series-mainthm}, every continuous endomorphism $\phi\in\End(\widehat{V}_{\mathbf{e}})$ is equal to the sum of the summable family $(\widehat{R}^{n}_{\mathbf{e},\boldsymbol{\alpha}}P_n(\widehat{L}_{\mathbf{e},\boldsymbol{\beta}}))_{n\in \N}$ for some polynomials sequence $(P_n)_{n\in\N}\in\K[\x]^{\N}$.
\end{corollaire}
\begin{proof}
Apply Theorem~\ref{series-mainthm} with $\boldsymbol{\beta}:=\boldsymbol{\alpha}\downarrow$ and
$\boldsymbol{\alpha}:=\boldsymbol{\beta}\uparrow$.
\end{proof}
\begin{note}
Without difficulty we can check that the extension $\widehat{D}_{\mathbf{e},\boldsymbol{\beta},\boldsymbol{\alpha}}$ of the diagonal operator $D_{\mathbf{e},\boldsymbol{\beta},\boldsymbol{\alpha}} = [L_{\mathbf{e},\boldsymbol{\beta}},R_{\mathbf{e},\boldsymbol{\alpha}}]$ is equal to $[\widehat{L}_{\mathbf{e},\boldsymbol{\beta}},\widehat{R}_{\mathbf{e},\boldsymbol{\alpha}}]$. As a continuous endomorphism,
$\widehat{D}_{\mathbf{e},\boldsymbol{\beta},\boldsymbol{\alpha}}=\displaystyle\sum_{n\in\N}\widehat{R}^{n}_{\mathbf{e},\boldsymbol{\alpha}}P_{n}(\widehat{L}_{\mathbf{e},\boldsymbol{\beta}})$. So the commutation rule becomes
\begin{equation}
\widehat{L}_{\mathbf{e},\boldsymbol{\beta}}\widehat{R}_{\mathbf{e},\boldsymbol{\alpha}}=\widehat{R}_{\mathbf{e},\boldsymbol{\alpha}}\widehat{L}_{\mathbf{e},\boldsymbol{\beta}}+\displaystyle \sum_{n\in\N}\widehat{R}^{n}_{\mathbf{e},\boldsymbol{\alpha}}P_{n}(\widehat{L}_{\mathbf{e},\boldsymbol{\beta}})\ .
\end{equation}
\end{note}

\section{Conclusions}
The idea of the commutation relation $AB-BA=I$ between two operators $A$ and $B$ (for example the creation and annihilation operators of
second-quantized theory)  is fundamental to the foundations of quantum physics.  In this paper we have shown that starting from this basic equality,
calculations of elementary operations, such as exponentiation associated with quantum dynamics and thermodynamics, lead us immediately to traditional
combinatorial concepts such as Stirling numbers, and generalizations thereof, which we describe.   We give explicit forms for the one-parameter groups
generated by the exponentials of such operators - crucial in quantum calculations -  in certain restricted cases; namely,  those  containing one-annihilator
only (corresponding to forms of so-called Sheffer-type).

In Physics, the creation and annihilation operators act on spaces of numbers of particles, moving from one state to another and so are considered as a special form of {\em ladder operator}.  We generalize this concept also, by considering endomorphisms in linear spaces, which mathematically correspond to these ideas. In particular, we note that infinite-dimensional vector space seems to be a rather natural
setting to deal with ladder operators. Any integer-indexed basis may provide the setting in a rather obvious
way for generalized ladder operators that can be either lowering (annihilation) or raising (creation), and
without any particular commutation rule. We prove that given two ladder operators, one lowering, the
other one raising, associated with possibly distinct bases (with the same first "rank"), it is possible to
expand any linear endomorphism in terms of iterates of the given ladder operators.

\end{document}